\documentclass[11pt]{amsart}

\usepackage{xy, amsfonts, amssymb,enumerate}

\newif\ifdraft\draftfalse
\newcommand\labl[1]{
  \ifdraft
     \ifmmode\renewcommand{\theequation}{\arabic{equation}-#1}
     \else{\bf(#1)}
     \fi
  \fi
  \label{#1}
}
\newcommand\Aff{\mathbb{A}}
 \newcommand\CC{\mathbb{C}}
  
 \newcommand\GG{\mathbb{G}}

\newcommand\PP{\mathbb{P}} \newcommand\QQ{\mathbb{Q}} 
\newcommand\RR{\mathbb{R}}

\newcommand\ZZ{\mathbb{Z}} 
\newcommand\Spec{\operatorname{Spec}}

\newtheorem{theorem}{Theorem}[section]
\newtheorem{lemma}[theorem]{Lemma}
\newtheorem{proposition}[theorem]{Proposition}
\newtheorem{corollary}[theorem]{Corollary}
\newtheorem{definition-proposition}[theorem]{Definition-Proposition}

\theoremstyle{definition}

\newtheorem{definition}[theorem]{Definition}
\newtheorem{example}[theorem]{Example}

\newtheorem{claim}[theorem]{Claim}

\theoremstyle{remark}
\newtheorem{remark}[theorem]{Remark}

\numberwithin{equation}{section}

\draftfalse \ifdraft \usepackage{showkeys} \fi
\usepackage[pagebackref=false]{hyperref}  

\newtheorem{exercise}[theorem]{Exercise}

\theoremstyle{remark}
\newtheorem{remarks}[theorem]{Remarks}

\newcommand\Chow{\operatorname{Chow}}

\newcommand\End{\operatorname{End}}
\newcommand\Corr{\operatorname{Corr}}

\newcommand\Hom{\operatorname{Hom}}
\newcommand\id{\operatorname{id}}

\newcommand\Char{{\operatorname{char\ }}}
\newcommand\Ad{{\operatorname{Ad}}}
\newcommand\PSO{{\operatorname{PSO}}}
\newcommand\sep{{\rm sep}}
\newcommand\ksep{k_{\sep}}
\newcommand\Gal{{\operatorname{Gal}}}
\newcommand\inn{{\rm inn}}
\newcommand\ad{{\rm ad}}
\newcommand\Ginn{G_{\rm inn}}
\newcommand\Gad{G_{\rm ad}}
\newcommand\pt{{\rm pt}}
\newcommand\pr{{\rm pr}}
\newcommand\GSch{$G$-{\rm Sch}}
\newcommand\Sch{{\rm Sch}}
\newcommand\QH{{\mathcal Q\mathcal H}}
\newcommand\nbar{\overline{n}}

\newcommand\IversenHTML{http://134.76.163.65/servlet/digbib?template=view.html&id=167779&startpage=235&endpage=242&image-path=http://134.76.176.141/cgi-bin/letgifsfly.cgi&image-subpath=/4257&image-subpath=4257&pagenumber=235&imageset-id=4257}

\begin{document}
\title[Motivic decompositions]{On motivic decompositions arising from the method of Bia{\l}ynicki-Birula}
\author{Patrick Brosnan}
\address{Dept. of Mathematics\\
North Campus\\
SUNY Buffalo\\
Buffalo NY 14260\\
USA}
\curraddr{School of Mathematics\\ Institute for Advanced Study\\
1 Einstein Drive\\Princeton NJ 08540\\ USA} 
\email{pbrosnan@buffalo.edu}
\begin{abstract}
  Recently, V.~Chernousov, S.~Gille and A.~Merkurjev have obtained a
  decomposition of the motive of an isotropic smooth projective
  homogeneous variety analogous to the Bruhat decomposition.  Using
  the method of A.~Bia{l}ynicki-Birula and a corollary, which is 
  essentially due to S.~del Ba\~{n}o, I generalize this decomposition
  to the case of a (possibly anisotropic) smooth projective variety
  homogeneous under the action of an isotropic reductive group.  This
  answers a question of N.~Karpenko.
\end{abstract}
\maketitle
\bibliographystyle{plain}  

\section{Introduction}  
An important difference between the category of motives and the
category of algebraic varieties over a field is the existence of interesting
direct sum decompositions of motives.  The simplest of these is the
decomposition of the Chow motive $M(\PP^n)$ of $n$-dimensional
projective space over a field $k$:
\begin{equation}
\label{eq:Pn}
M(\PP^n)=\oplus_{i=0}^n \ZZ(n).
\end{equation}
This is one of the elementary results in Grothendieck's theory of motives (which will be recalled 
in \S\ref{sec-motives}).   
 
An example of a less elementary decomposition is the following theorem
due to M.~Rost, which is in an important ingredient
in his construction of the ``Rost motive''~\cite[Proposition 2]{RostPfister}.
\begin{theorem}[Rost decomposition]
Let $Q$ be a smooth, projective, $n$-dimensional isotropic quadric
over a field $k$ of characteristic not equal to $2$.  Then
\begin{equation}
\label{eq:Quadric}
M(Q)=\ZZ\oplus M(Q')(1)\oplus \ZZ(n)
\end{equation}
where $Q'$ is a smooth, projective sub-quadric of codimension $2$ in
$Q$.
\end{theorem}

Since both projective spaces and quadrics are examples of projective
homogeneous varieties, it is natural to look for decompositions
generalizing~(\ref{eq:Pn}) and~(\ref{eq:Quadric}) in the motives of
such varieties.  In the case that $G$ is a split reductive group (i.e.,
when the base field $k$ is separably closed), a decomposition
for the motive of $G/P$ was found by B.~K\"{o}ck~\cite{Kock}.  In this
case $M(G/P)$ splits completely as a sum of Tate motives.  A more
general decomposition was later found by
N.~Karpenko~\cite{karpenko-flag} in the case of motives of flag
varieties for classical groups, and recently K\"{o}ck's decomposition was
further generalized by V.~Chernousov, S.~Gille and
A.~Merkurjev~\cite[Theorem 7.4]{CGM2003} to the case of motives for
isotropic projective homogeneous varieties for adjoint semi-simple
groups.  Both generalizations explicitly describe the factors
appearing in the decomposition (which are, in general, not Tate
motives) in terms of smaller projective homogeneous spaces.  Moreover,
Rost's theorem appears as a special case of either generalization when
the quadric $Q$ is viewed as a homogeneous space for the group
$\PSO(q)$ with $q$ a quadratic form whose corresponding projective
quadric is $Q$.

Theorem~\ref{thm:main} of this paper applies results of
Bia{\l}ynicki-Birula (as extended by W.~Hesselink) to obtain a
decomposition of the motive of any smooth, projective algebraic
variety admitting an action of the multiplicative group.  In the case
of motives with rational coefficients (see Remark~\ref{rem-distinction}
for the distinction), the theorem is due to del
Ba\~{n}o~\cite{delBano}.  (I thank B.~K\"{o}ck for bringing this to my
attention.)  I include a sketch of the proof (which is similar to del
Ba\~{n}o's) for completeness and convenience of the reader.

As the class of varieties admitting multiplicative group actions
includes both the homogeneous spaces considered by Karpenko and the
isotropic homogeneous spaces considered by Chernousov, Gille and
Merkurjev, we obtain a generalization of the
Chernousov-Gille-Merkurjev decomposition.  In the end, we can give a
rather explicit decomposition of the motive $M(X)$ of a projective
homogeneous variety for a reductive group $G$ as a sum of Tate twisted
motives of certain ``quasi-homogeneous schemes'' $X_i$  for the
anisotropic kernel of $G$.  
Theorem~\ref{thm:RoughAnswer} gives a
rough form of this description which is refined in
Theorem~\ref{thm-Final}, the final theorem of the paper.  In
particular, the fact that the motive of a projective homogeneous
variety admits a decomposition in terms of motives of
quasi-homogeneous schemes for the anisotropic kernel answers the
fundamental question posed by Karpenko in the introduction
to~\cite{karpenko-flag}.

\subsection{Notation}
  All notions of Chow groups are taken from Fulton's book on 
  intersection theory~\cite{fulton}.  The official reference for
  reductive groups is SGA3~\cite{SGA3.3}, but some notation is taken
  from Springer's book~\cite{Springer}.  The main difference between
  these two references that will be important for this paper is that
  SGA3 demands that a reductive group is {\em connected} and a
  parabolic subgroup is {\em smooth}, while Springer does not make
  these assumptions.  Since these are convenient assumptions for us,
  we will have to agree with SGA3.  In several places, the symbol
  ``$k$'' is used to denote both the base field and an index.  This
  does not seem to produce any confusion.

\subsection{Outline}  With one exception, the results in this 
paper build sequentially from general facts about motivic
decompositions to specific information about the motivic
decomposition of a projective homogeneous space in terms of double
cosets of the Weyl group given in Theorem~\ref{thm-Final}.
Specifically, \S\ref{sec-motives} reviews the theory of motives,
\S\ref{sec-motivic-decomposition} explains how Bia{\l}ynicki-Birula's
theorem yields a motivic decomposition and
\S\ref{sec-phsc} introduces the concept of projective quasi-homogeneous
schemes and studies motivic decompositions of such schemes in the
presence of a $\GG_m$-action.  The last two mathematical sections,
\S\ref{sec-weyl} and \S\ref{sec-explicit}, formulate the general theory
in terms of root systems and reflection groups.  The one exceptional
section is \S\ref{sec-nilpotence} where I give a proof of a
generalization of Rost's nilpotence theorem following Chernousov,
Gille and Merkurjev.  The results in this section are not needed
anywhere else in the paper.  

\section{Motives and the category of correspondences}
\label{sec-motives}
The category of motives can be defined by first defining the category of
correspondences and then applying the functor of idempotent
completion.  In fact, the decomposition theorems of this paper such as
Theorem~\ref{thm:main} will hold before (or after) taking idempotent
completion, but, to make the connection with Chow motives explicit, I
will describe both categories.

The category of correspondences is the category $\Corr_k$ whose objects
are pairs $(X,n)$ with $X$ a smooth projective scheme over the field
$k$ and $n$ an integer.  The morphisms are given by
\begin{equation}
\label{eq:chowmotives}
Hom_{\Corr_k}((X,n), (Y,m))= \oplus A_{d_i + n-m} X_i\times Y
\end{equation} 
where $X=\coprod X_i$ with $X_i$ connected, $d_i=\dim X_i$ and $A_k$ denotes the $k$-th Chow groups graded by dimension.

For a smooth, projective, variety $X$, $M(X)$ denotes the object
$(X,0)$.  When $M=(X,n)$ is an object, $M(k)$ denotes the {\it Tate
twisted} object $(X,n+k)$.  Let $\ZZ$ denote the object $(\Spec k, 0)$. The ``twists'' $\ZZ(k)$
of $\ZZ$ are called {\it Tate objects}.    Clearly, $\Hom_{\Corr_k} (Z(k), M(X))= A_kX$ and 
$\Hom_{\Corr_k} (M(X), Z(k)) = A_{d-k}X$ for $X$ irreducible of dimension $d$.

The objects of the category $\Chow_k$ of Chow motives over $k$ are
triples $(X,n,p)$ with $p\in\End (X,n)$ a morphism such that $p^2=p$.
The morphisms in $\Chow_k$ are given by
$$
Hom_{\Chow_k}((X,n,p),(Y,m,q))= qHom_{\Corr_k} ((X,n),(Y,m))p.
$$ The category $\Chow_k$ is {\it idempotent complete} in the sense
that every idempotent morphism has both a kernel and a cokernel.
Clearly, there is a fully, faithful embedding
$\Corr_k\hookrightarrow\Chow_k$ given by $(X,n)\rightsquigarrow
(X,n,\id)$.  Moreover, this embedding is universal for functors from
$\Corr_k$ with idempotent complete targets.

Both categories admit a tensor structure defined (on $\Chow_k$) by
$$
(X,n,p)\otimes (Y,m,q)=(X\times Y, n+m, p\times q).
$$

The category $\Chow_k$ also admits direct sums defined as follows.   Let $r_k$ denote the idempotent on 
$M(\PP^k)$ such that $\ZZ(k)=(\PP^k,0,r_k)$.   Explicitly, it is given by the cycle $[\pt \times \PP^k]$
 where $\pt$ denotes an arbitrary degree $1$ closed point in $\PP^k$.  Then, for two motives $(X,n,p)$ and $(Y,m,q)$ with $n\leq m$,
$$
(X,n,p)\oplus (Y,m,q)=(X\coprod (Y\times\PP^{m-n}), n, p +  (q\times r_{m-n})).
$$ 
The direct sum is the coproduct in the category of motives.  In the
category of correspondences, coproducts do not always exist.  For
example, it is not hard to see that the object $\ZZ\coprod \ZZ(2)$
does not exist in the category of correspondences over $\CC$.
However, when coproducts do exist in $\Corr_k$, they coincide with
those in the category of Chow motives.

\begin{remark}
\label{rem-distinction}
The categories of motives occurring in this paper (and in those papers
of Rost, Karpenko, and Chernousov-Gille-Merkurjev) have {\it integral
  coefficients}.  If we were to tensor the morphism sets with $\QQ$,
replacing $\Hom_{\Corr_k} (M(X),M(Y))$ with
$\Hom_{\Corr_k}(M(X),M(Y))\otimes\QQ$, we would obtain a category $\Corr_k\otimes
\QQ$ which is closer to the categories of motives Grothendieck
originally considered~\cite{Jannsen, Manin}.  However, we would also
loose information.  For example, using the fact that any quadric is
totally isotropic over a finite separable extension, it is easy to see
that every quadric decomposes as a direct sum of the Tate objects
$\QQ(i)=\ZZ(i)\otimes \QQ$ in $\Corr_k\otimes\QQ$.  On the other hand, Springer's
theorem on quadrics isotropic over an odd degree
extension~\cite[Theorem 2.3 p. 198]{Lam} implies that the integral
motive $M(Q)\in \Corr_k$ of a smooth quadric contains a factor of
$\ZZ(0)$ if an only if $Q(k)\neq \emptyset$.
\end{remark}

\section{Motivic decomposition}
\label{sec-motivic-decomposition}

The most general theorem in this paper on motivic decompositions is essentially
a corollary of two results which I will now recall after giving one
definition.   As mentioned in the introduction, a version of the theorem (\ref{thm:main}) can also be found
in Theorem~2.4 of S.~del Ba\~{n}o's paper, ~\cite{delBano}. 

\begin{definition}
A flat morphism $\phi:X\to Z$ is called an {\it affine fibration} (resp. an {\it affine quasi-fibration}) of
relative dimension $d$ if,
for every point $z\in Z$, there is a Zariski open neighborhood $U\subset Z$ such that $X_U\cong Z\times\Aff^d$
with $\phi:X_U\to Z$ isomorphic to the projection on the first factor (resp. the fiber $X_z$
of $\phi$ is isomorphic to $\Aff^d_{k(z)}$).  
\end{definition} 

Clearly an affine fibration is an affine quasi-fibration.  
It is a well-known consequence of the homotopy invariance of Chow
groups that an affine quasi-fibration between smooth varieties of
relative dimension $d$ induces an isomorphism $\phi^*:A_i Z\to A_{i+d}
X$.  

\begin{theorem}[Karpenko] 
\label{thm-karpenko}
Let $X$ be a smooth, projective variety over a field $k$ with a filtration
$$ X=X_n\supset X_{n-1}\supset\dots\supset X_0\supset X_{-1}=\emptyset
$$
where the $X_i$ are closed subvarieties.
Assume that, for each integer $i\in [0,n]$, there is a smooth projective variety $Z_i$ and an
affine fibration 
$\phi_i:X_i-X_{i-1}\to Z_i$
of relative dimension $a_i$.   Then, in the category of correspondences, 
$M(X)=\coprod_{i=0}^n M(Z_i)(a_i)$.
\end{theorem}

The theorem was stated by Karpenko~\cite{karpenko-flag} for the special case that the maps
$\phi_i:X_i-X_{i-1}\to Z_i$ are vector bundle morphisms.  However, in~\cite[Theorem 7.1]{CGM2003},
Chernousov, Gille and Merkurjev noticed that Karpenko's 
proof actually applies to any affine quasi-fibration.  (Del Ba\~{n}o gives a slightly different proof of the result
in the proof of his~Theorem~2.4).

The second result is the method of Bia{\l}ynicki-Birula which gives a natural situation
where Karpenko's theorem applies.

\begin{theorem}[Bia{\l}ynicki-Birula, Hesselink, Iversen]
\label{thm-bb}
Let $X$ be a smooth, projective variety over a field $k$ equipped with
an action of the multiplicative group $\GG_m$.  Then
\begin{enumerate}
\item The fixed point locus $X^{\GG_m}$ is a smooth, closed subscheme of $X$.
\item There is a numbering $X^{\GG_m}=\coprod_{i=1}^n Z_i$ of the connected components of the fixed point
locus, a filtration 
$$ X=X_n\supset X_{n-1}\supset\dots\supset X_0\supset X_{-1}=\emptyset
$$
and affine fibrations $\phi_i:X_i-X_{i-1}\to Z_i$.   
\item The relative dimension $a_i$ of the affine fibration $\phi_i$ is the dimension
of the positive eigenspace of the action of $\GG_m$ on the tangent space of $X$ at an arbitrary
point $z\in Z_i$.  The dimension of $Z_i$ is the dimension of $TX_z^{\GG_m}$.
\end{enumerate} 
\end{theorem}

As the theorem stated is the product of several results of different
authors, I will give a short history of the result (in lieu of a proof).

\href{\IversenHTML}{Iversen}~\cite{IversenFixed} showed that $X^{\GG_m}$ is smooth.  
Bia{\l}ynicki-Birula~\cite{BialynickiBirulaDecomposition}
showed, under the assumption that $k$ is algebraically closed, that
$X$ is a union of locally closed subschemes $X_i^{+}$ with affine
fibrations $\phi_i:X_i^+\to Z_i$ where the $Z_i$ are the connected
components of $X^{\GG_m}$.  In fact, this was shown with the assumption that $X$ is 
projective replaced with the assumption that $X$ is complete.
Shortly thereafter, Bia{\l}ynicki-Birula showed
that, when $X$ is projective, there is a filtration of $X$ and an
ordering of the connected components as in the theorem such that
$X_i^+=X_i-X_{i-1}$~\cite{BialynickiBirulaFiltration}.  Thus, in the case that $k$ is algebraically
closed, the theorem as stated was proved by Bia{\l}ynicki-Birula and
Iversen. Note that the existence of a filtration is deduced by
embedding $X$ equivariantly in a projective space with a
diagonalized $\GG_m$-action.  It is then easy to construct a
filtration on the projective space and see that that it restricts to
one on $X$.    (However, for $X$ smooth and complete but not projective, there are
examples where no filtration satisfying $X_i^+=X_i-X_{i-1}$ exists~\cite[Example 2, p. 30]{CarrellSommese}.)

Bia{\l}ynicki-Birula's proofs make use of the assumption that $k$ is
algebraically closed.  Hesselink removed this restriction and was able
to show that $X$ is a union of locally closed $X_i^{+}$ for $X$ a
smooth, proper scheme over an arbitrary base~\cite{HesselinkConcentration} provided that 
there is a covering of $X$ by $\GG_m$-stable Zariski open affine subsets.  To construct a
filtration in the case $X$ is projective, it then suffices to find a
$\GG_m$-equivariant embedding of $X$ into a projective space with a
diagonal action of $\GG_m$ or, equivalently, a very ample
$\GG_m$-linearized line bundle over $X$.  The fact that such bundles
exist can be found in Mumford's GIT~\cite{GIT}.   Since a diagonal action of $\GG_m$ on $\PP^r$ 
preserves the coordinate hyperplanes (and thus their complements), any $X$ embedded by a $\GG_m$-linearized
very ample line bundle automatically has a covering by $\GG_m$-invariant open affines.  Thus Hesselink's 
hypotheses are verified for $X$ smooth and projective.

The Bia{\l}ynicki-Birula decomposition is explicit in the sense that
the locally closed subscheme $X_i^+$ is the set of all points $x\in X$
such that $\lim_{t\to 0} tx\in Z_i$ where $(t,x)\mapsto tx$ is the
$\GG_m$ action.  Moreover, the map $\phi_i:X_i^+\to Z_i$ is then given
by $x\mapsto \lim_{t\to 0} tx$. 

\begin{remark}
Since $X$ separated, $\lim_{t\to 0} tx$ has at most one meaning, since $X$ is 
proper it has exactly one.
\end{remark}

\begin{theorem}
\label{thm:main}
Let $X$ be a smooth, projective scheme over a field $k$ equipped with
an action of the multiplicative group $\GG_m$.  Then, in the category
$\Corr_k$,
\begin{equation}
\label{eq:main}
M(X)=\coprod M(Z_i)(a_i)
\end{equation}
where the $Z_i$ are the connected components of $X^{\GG_m}$ and the
$a_i$ are determined as in Theorem~\ref{thm-bb} (3).
\end{theorem}
\begin{proof}
This is a corollary of the two previous theorems.
\end{proof}

We will refer to the decomposition of \eqref{eq:main} as the {\em motivic Bia{\l}ynicki-Birula decomposition.} 
The rest of this paper will focus on the application of the theorem to
the special case where $X$ is a projective homogeneous variety.
Before proceeding to the general theory, I use the theorem directly to
derive Rost's decomposition theorem.

\begin{example}
\label{ex:IsotropicQuadric}
 Let $q:V\to k$ be a non-degenerate quadratic form of 
  dimension $n+2$ over a field $k$ with $\Char k\neq 2$, and let $Q$ be the associated
$n$-dimensional smooth projective quadric.    Suppose that
  $q$ is isotropic, that is, there exists a nonzero vector $v\in V$
  such that $q(v)=0$.  Then there is a subspace $W\subset V$ and two
  linearly independent vectors $v_1$ and $v_2$ such that
  $V=kv_1\oplus kv_2\oplus W$ and
$$
q(xv_1 + yv_2 + w) = xy+q'(w)
$$
for a non-degenerate quadratic form $q'$ on $W$.  
(This is an
easy exercise which can also be found in almost any book on quadratic
forms over a field, e.g., \cite[Proposition 3.7.1]{Knus}.)
In this case, the
multiplicative group $\GG_m$ acts on $Q$ by
$$
[xv_1+yv_2+w]\mapsto [txv_1+t^{-1}yv_2+w].
$$
Assume that $\dim Q\neq 2$.  
The fixed point set $Q^{\GG_m}$ then has three components: the points $[v_1]$ and $[v_2]$ and the 
quadric $Q'=\{[w]\in\PP(W)\,|\, q'(w)=0\}$ which we denote by $Z_1$, $Z_2$ and $Z_3$ 
respectively.  
Let $T_i$ denote the tangent space at an (arbitrary) point of $Z_i$.  
The action of $\GG_m$ on $T_1$ has only negative weights.  Therefore, in the decomposition of  
Theorem~\ref{thm:main}, $a_1=0$.  The weights of $T_2$ are all positive, therefore, $a_2=n$.  
Finally, $T_3$ has weights $-1$ and $1$ each occurring once and $0$ occurring $n-2$ times.  Therefore
$a_3=1$.  Thus we have 
\begin{equation}
\label{eq:Rost2}
M(Q)=\ZZ\coprod\ZZ(n)\coprod M(Q')(1).
\end{equation}
It is easy to see that the above decomposition also holds when $\dim
Q=2$, however there is a possibility that the quadric $Q'$ can split
as a disjoint union of $2$ copies of $\Spec k$.
\end{example}

\section{Projective homogeneous schemes}
\label{sec-phsc}

Let $G$ denote a {\em reductive} group over a field $k$ in the terminology of SGA3.   That is, $G$ is smooth and connected 
with trivial unipotent radical.   Recall that a {\em parabolic} subgroup of $G$  is a subgroup $P$ such that $G/P$ is projective and 
$P$ is smooth over $k$.   Let $\GSch_k$ denote the category of  $G$-schemes over $k$.   The objects of this category
are schemes $X$ over $k$ equipped with a $G$-action $G\times X\to X$.   The morphisms are the $G$-equivariant 
scheme-theoretic morphisms.  Base change induces an obvious functor $\GSch_k\leadsto \GSch_L$ for $L$ an extension
of $k$.

Let $\overline{k}$ denote an algebraic closure of $k$.  A $G$-scheme
$X$ is a {\em projective homogeneous variety} for $G$ if
$X_{\overline{k}}$ is isomorphic as a $G_{\overline{k}}$-scheme to
$G_{\overline{k}}/P$ for $P\subset G_{\overline{k}}$ a parabolic subgroup.
It is well-known that such a projective homogeneous variety is
projective over $k$.  We will call a $G$-scheme $X$ a {\em projective
  quasi-homogeneous scheme} if $X$ is smooth and projective over
$k$ and the morphism $\psi=(a,\pr_2):G\times X\to X\times X$ given by
$(g,x)\mapsto (gx,x)$ is smooth.

\begin{proposition}\label{prop:quasi-hom}   Let $X$ be a $G$-scheme over $k$.   Then the following are equivalent.
  \begin{enumerate}
  \item $X$ is a projective quasi-homogeneous $G$-scheme.
  \item $X_{\overline{k}}$ is a disjoint union of projective homogeneous varieties.
  \item $X$ is smooth, projective, and, for every geometric point $x\in
    X_{\overline{k}}$, the orbit map $m_x: G_{\overline{k}}\to
    X_{\overline{k}}$ (given by $g\mapsto gx$) induces a surjection $dm_x:
    L(G)\to TX_x$.
  \end{enumerate}
\end{proposition}
\begin{proof}
Since all of the  properties listed are invariant under base change of $k$, we can assume that $k$ is algebraically closed.  

(1) $\Rightarrow$ (2):  A scheme is projective
quasi-homogeneous if and only if all of its connected components are
projective quasi-homogeneous.   (This is easy.)  Therefore we can assume that $X$ is
connected.  Since $G\times X\to X\times X$ is smooth, all orbits are open.  Thus,
since $X$ is smooth and connected (hence irreducible), all orbits must
intersect.  It follows that there is only one orbit, namely, $X$
itself.  Thus $X=G/P$ for some subgroup $P$.  The smoothness of $\psi$
then implies that $P$ is smooth.

(2) $\Rightarrow$ (3):  Here we can assume $X=G/P$ with $P$ parabolic.   The claim then follows from the assumption that 
$P$ is smooth.

(3)$\Rightarrow$ (1):  From (3), it follows  that
$\psi:G\times X\to X\times X$ induces surjections on the tangent spaces.  The claim then
follows from \cite[p. 270, Proposition 10.4 (iii)]{HartshorneAG}.
\end{proof}

\begin{remark}
  The motivation for considering projective quasi-homogeneous schemes in addition to projective
homogeneous varieties is already apparent in the Rost decomposition of quadrics.    Let $q$ be an 
$n$-dimensional non-degenerate quadratic form, let $Q$ be the associated $n-2$ dimensional 
smooth projective quadric and let $\PSO(q)$ denote the special orthogonal group.   Then 
$Q$ is a projective homogeneous space for $\PSO(q)$ if and only if $\dim Q>0$.     When $\dim Q=0$,
$Q$ can either be irreducible or a disjoint union of two copies of $\Spec k$.    In either case,
it is not a projective homogeneous variety for $\PSO(q)$.    However, regardless of the dimension, 
$Q$ is projective quasi-homogeneous for $\PSO(q)$.
\end{remark}

Now let $X$ denote a projective quasi-homogeneous scheme for $G$, and let $\GG_m\stackrel{\lambda}{\stackrel{\sim}{\to}} L\subset G$
denote the inclusion of a $k$-split torus in $G$.   (The group $G$ is called {\em isotropic} if 
such a split torus exists.)
In this case, $L$ acts on $X$ and Theorem~\ref{thm:main} applies to give
a decomposition of $M(X)$.  The main result of this section is that, in fact, the summands
appearing are themselves projective quasi-homogeneous schemes for the centralizer $H=Z(\lambda)$ of $L$ in $G$.
A more detailed description will be obtained 
in~\S\ref{sec-weyl} and \S\ref{sec-explicit}.

\begin{theorem}
\label{thm:description}
  \begin{enumerate}
  \item $H$ is connected, reductive and defined over $k$. 
  \item $H$ acts on the fixed point set $X^{\lambda}$.
  \item The action map $\psi_H: H\times X^{\lambda} \to X^{\lambda}\times X^{\lambda}$ is smooth.   Thus $X^{\lambda}$ is 
projective quasi-homogeneous.
  \end{enumerate}
\end{theorem}
\begin{proof}
  (1)  The fact that $H$ is connected is [Springer, 13.4.2 (i)].  The
  fact that it is reductive is [Springer, 7.6.4 (i)].  It is defined
  over $k$ by [Springer, 13.3.1 (ii)].
  
  (2) To see that $H$ acts on $X^{\lambda}$, let $T$ be a scheme over $k$ and consider $T$-valued points $x\in X^{\lambda}(T)$, $h\in H(T)$
  and $t\in L(T)$.  Then $thx=htx=hx$, thus, $hx$ is in $X^{\lambda}(T)$.

(3) Since the smoothness of $\psi_H$ is invariant under field extension of $k$, we may assume that $k=\overline{k}$.
Pick a closed point $z\in Z$.   The orbit map $l:G\to X$ given by $g\mapsto gz$, induces a surjection
$dl:L(G)\to TX_z$ because, by assumption, $X$ is a projective quasi-homogeneous scheme for $G$.  
The multiplicative group $\GG_m$ acts on $G$ via conjugation by $\lambda$, i.e., $g\mapsto \lambda(t)g\lambda(t^{-1})$.
The group $\GG_m$ also acts on $X$ via right multiplication, i.e., $x\mapsto \lambda(t)x$.
Since $z$ is a fixed point of $\lambda$, the orbit map $l$ is equivariant for the $\GG_m$-actions.   Moreover, we obtain 
a $\GG_m$-action on $TX_z$ compatible with the $Ad$-action of $\GG_m$ on $L(G)$.  

Now $L(G)\cong L(G)_+ \oplus L(H)$ where $L(G)_+$ consists of the non-zero weight space of $L(G)$ and $L(H)$ is the Lie
algebra of $H$.   Analogously, $TX_z\cong TX_{z+}\oplus TX^{\lambda}_z$ where $TX_{z+}$ is the non-zero weight space of $TX_z$.
Since $dl$ respects the weight decomposition, $dl(L(H))=TX^{\lambda}$.   Thus Proposition~\ref{prop:quasi-hom} (3)
is satisfied. 
\end{proof}

\begin{corollary}
\label{cor:description}
  Let $X$ be a projective quasi-homogeneous scheme for an isotropic reductive group $G$, and let $\lambda:\GG_m\to G$ be the embedding
of a split torus.  Then, in the motivic Bia{\l}ynicki-Birula decomposition
  \begin{equation*}
    M(X)=\coprod M(Z_i)(a_i)
  \end{equation*}
of Theorem~\ref{thm:main}, the $Z_i$ are all projective quasi-homogeneous schemes for the centralizer $H$ of $\lambda$.
\end{corollary}
\begin{proof}
The corollary holds because the $Z_i$ appearing in Theorem~\ref{thm:main} are components of $X^{\lambda}$.
\end{proof}

\subsection{Adjoint groups} For a reductive group $G$, let $\QH_G$ denote the full subcategory of $G-\Sch_k$
consisting of
projective quasi-homogeneous schemes.  If $Z_G$ is
the center of $G$, then $G_{\ad}=G/Z_G$ is the {\em adjoint} group of $G$,
an adjoint semi-simple group~\cite[Proposition 22.4.3.5]{SGA3.3}.  The
restriction functor $G_{\ad}-\Sch_k\rightsquigarrow G-\Sch_k$ induces
an equivalence of categories $\QH_{G_{\ad}}\rightsquigarrow \QH_G$.
This is because $Z_G$ is smooth and acts trivially on all
quasi-homogeneous schemes over $G$. 

\subsection{Anisotropic Kernels}\cite{TitsBoulder}
Let $S$ denote
a maximal $k$-split torus of $G$, and let $Z(S)$ denote its centralizer. The derived subgroup $DZ(S)$ is the {\em semi-simple
anisotropic kernel}. 
Since $Z(S)$ is reductive, there is an almost direct product
decomposition $Z(S)=DZ(S)\cdot Z$ where $Z$ is the center of
$Z(S)$\cite[Proposition 2.2]{borel-tits}.   It follows that the adjoint group of $Z(S)$ is isomorphic to the adjoint group of the 
semi-simple anisotropic kernel.   Thus the categories $\QH_{Z(S)}$, $\QH_{DZ(S)}$ and $QH_{Z(S)_{\ad}}$ are all equivalent 
via the restriction of group functors.   (Moreover, the objects are identical.) 

Applying Corollary~\ref{cor:description} inductively, we obtain an
answer to a question of N.~Karpenko~\cite{karpenko-flag}.

\begin{theorem}\label{thm:RoughAnswer}
  Let $X$ be a projective quasi-homogeneous scheme for a reductive
  group $G$.  Then
\begin{equation}
\label{eq:KarpenkoAnswer}
M(X)=\coprod M(Y_i)(a_i)
\end{equation}
where the $Y_i$ are 
irreducible projective quasi-homogeneous schemes for the anisotropic kernel of $G$ (resp. for $Z(S)$, for $Z(S)_{\ad}$).  
\end{theorem}

\begin{remark}
  A projective homogeneous variety $X$ is said to be {\em isotropic}
  if $X=G/P$ for a parabolic subgroup $P$ defined over $k$.  Otherwise it is said to be {\em anisotropic}.
  $X$ is anisotropic if and only if $X(k)$ is empty.
  If $X$ is
  an isotropic projective homogeneous space for a reductive group $G$,
  then there exists at least one $k$-split torus $L$ in $G$.
  (See~\cite{Springer} or~\cite{SGA3.3}.)   In other words, if $X$ is
  isotropic then $G$ is as well.  It follows that the schemes $Y_i$ appearing in (\ref{eq:KarpenkoAnswer}) are all
either anisotropic or isomorphic to $\Spec k$.
\end{remark}

\begin{exercise}
\label{rem-InterestingExample}
  It is interesting to see an example of 
  Corollary~\ref{cor:description} at work on an anisotropic
  projective homogeneous variety for an isotropic reductive group.  One
  such is given by the variety $X$ of two dimensional isotropic
  subspaces for the quadratic form $q=x_1^2+\cdots +x_{2n}^2 + yz$ over the
  reals with $n\geq 2$.   Using the methods of Example~\ref{ex:IsotropicQuadric}, the 
  decomposition can be computed explicitly in terms of smaller quadrics and varieties
  of two dimensional isotropic subspaces.
In Example~\ref{example-TwoDIsotropic}, we will return to this matter, computing the decomposition 
using the Lie theory of $\PSO(q)$.
\end{exercise}

\section{The nilpotence theorem of Chernousov, Gille and Merkurjev}
\label{sec-nilpotence}

As a corollary of the results of the previous section, we obtain the following theorem of
Chernousov, Gille and Merkurjev~\cite[Theorem 8.2]{CGM2003}.
\begin{theorem}
  Let $X$ be a projective homogeneous variety for a reductive group $G$
  over a field $k$.  Then the kernel of the map
$$
\End(M(X))\to \End(M(X\otimes\overline{k}))
$$
consists of nilpotent endomorphisms.
\end{theorem}

The proof follows that of~\cite{CGM2003}.  I include it here
for the convenience of the reader and to make the point that the
theorem can be obtained without the full description of the motivic
decomposition obtained in \cite[Theorem 7.4]{CGM2003}.

For a field extension $L/k$, let $n_L$ denote the number of terms appearing in the decomposition
\begin{equation}
\label{eq:decomposition}
M(X_L)=\coprod_{i=1}^{n_L} M(Z_i)(a_i)
\end{equation}
of (\ref{eq:KarpenkoAnswer}) for the projective homogeneous $G_L$-variety $X_L$.  (Here the $Z_i$ will
depend on $L$.)
Clearly, $M\supset L\Rightarrow n_M\geq n_L$, and the maximal
number of terms in the coproduct occurs precisely when each $Z_i$ is
$\Spec L$.  In particular, this happens when $L=\Spec \ksep$.

\begin{claim}
  Set $N(d,n)=(d+1)^{\nbar-n}$ with $\nbar=n_{\ksep}$.
  Then, for any morphism
  $f\in\End(M(X))$ with $f\otimes\overline{k}=0$, $f^{N(d,n_k)}=0$.  
\end{claim}  
 
  Evidently the claim implies the theorem.

  Now in the case that the maximal number of terms appears in the
  decomposition (i.e., $n_k=\nbar$), the claim is trivial because each of the objects
  appearing is Tate.  In fact, the only
  morphism in $\End(M(X))$ which vanishes in
  $\End(M(X)\otimes\overline{k})$ is the $0$ morphism.  Thus the claim is valid for $n_k=\nbar$.
  
  Now reason by descending induction on $n=n_k$.  (Properly speaking,
  we reason by ascending induction on $\nbar -n_k$ starting with the
  case $\nbar-n_k=0$.)  Let $f\in\End(M(X))$ be an endomorphism in the kernel 
of the map to $\End(M(X\otimes\overline{k}))$ and pick
  a point $z$  in one of the anisotropic
  components $Z_i$ appearing in (\ref{eq:decomposition}).  (If all
  components are isotropic, $n$ is maximal and the claim is already
  proved.)  Set $L=k(z)$.  Over $L$, $Z_i$ is isotropic.  Therefore
  the number $n_i=n_L$ of terms appearing in the motivic decomposition
  of $X_L$ is greater than $n$.  Thus the claim holds for $X_L$ and
  $f_L^{N(d,n_i)}=0$.   Since $N(d,n_i)\leq N(d,n+1)$, it follows that $f_L^N=0$
where $N=N(d,n+1)$.
  
  I now use~\cite[Theorem 3.1 and Remark
  3.2]{brosnan-documenta} in the form in which it was used by Gille,
  Chernousov and Merkurjev~\cite[Proposition 8.1]{CGM2003}.

  \begin{lemma}
    Let $X$ be a smooth, projective variety over a field $k$ and $Z$
    an $r$-dimensional scheme of finite type over $k$.  Let
    $f\in\End(M(X))$ be an endomorphism such that, for every point
    $z\in Z$, the morphism $f_{z*}:A_*(X\otimes k(z))\to A_*(X\otimes
    k(z))$ vanishes.  Then $f_*^{r+1}:A_*(X\times Z)\to A_*(X\times
    Z)$ vanishes.
  \end{lemma}

From the lemma and the fact that $\dim Z_i\leq d$, it follows that the composition
\begin{equation}
  \label{eq:composition}
  M(Z_i)(a_i)\stackrel{j_i}{\to} M(X)\stackrel{f^{(d+1)N(d,n+1)}}{\to} M(X)
\end{equation}
vanishes where the first arrow, $j_i$,  is the canonical one coming from the
coproduct decomposition.  Thus, for each anisotropic $Z_i$ in the coproduct,
$f^{(d+1)N}\circ j_i=0$.  On the other hand, if $Z_i=\Spec k$, then it is easy
to see that the composition in \eqref{eq:composition} vanishes even
with $f^{(d+1)N(d,n+1)}$ replaced by $f$.    

Since $N(d,n+1)=(d+1)^{\nbar-n_k-1}$, the claim is
proved.

\begin{remark}
Clearly the exponent $(d+1)^{\nbar-n_k}$ is not optimal.  
\end{remark}

\section{The Weyl group and its double cosets}
\label{sec-weyl}

In this section and the next, I give an explicit description of the components
$Z_i$ appearing in the motivic decomposition \eqref{eq:main} of a
projective homogeneous variety $X$ for an isotropic reductive group
$G$.  Roughly speaking, the geometric components are in one-to-one
correspondence with certain double cosets of the Weyl group.  The
algebraic components correspond to equivalences classes of these
double cosets under the so-called ``*-action'' of the absolute Galois group of the
base field $k$~\eqref{subsec-StarAction}.

While the language of schemes was used in the previous sections, in
this section I abuse notation slightly (e.g. in the proof of
lemma~\ref{lemma:parabolic-normalizer}) and confuse points with
$k$-valued points.  This facilitates comparison with the
reference~\cite{Springer} which is written in the language of
varieties.

For $T\subset G$ a maximal (but not necessarily split) torus
defined over $k$, set $W=W(G,T)$, the corresponding Weyl group.  For a
subtorus $C\subset T$, let $W_C=W(Z_G(C),T)$ where $Z_G(C)$ denotes the
centralizer of $C$.  It is a subgroup of $W$.  Likewise, for a
character $\phi:\GG_m\to T$, let $W_{\phi}=W(Z_G(\phi), T)$,
also a subgroup of $W$.

Now, in the situation of the \S\ref{sec-phsc}, $G$ has a cocharacter
$\lambda:\GG_m\stackrel{\cong}{\to}L\subset G$.  We can assume that $L\subset S\subset T$ with 
$S$ a maximal $k$-split torus and $T$ a maximal torus with $T$ defined over 
$k$~\cite[Theorem 13.3.6 and Remark 13.3.7]{Springer}.

If $X$ is isotropic, there is a parabolic subgroup $P$ such that
$X=G/P$.  Since $G$ is reductive, we may assume that $P=P(\mu)$ for a
cocharacter $\mu:\GG_m\to S$.   (Roughly, $P(\mu)$ is defined as the set of $g\in G$ such that 
$\lim_{t\to\infty} \mu(t)g\mu(t)^{-1}$ exits.  See~\cite[\S13.4.1]{Springer}.)  Set
\begin{equation}
  \label{eq:pi_not}
\mathcal{X}=  W_{\lambda}\backslash W/W_{\mu}.
\end{equation}

\begin{theorem}
\label{thm:DoubleCosets}
  If $k$ is separably closed, the connected components $Z_i$ appearing
  in the motivic Bia{\l}ynicki-Birula decomposition of $X$ are in one-one correspondence
  with the elements of $\mathcal{X}$.
\end{theorem}

To begin the proof of the theorem, first note that, since $k$ is
separably closed, $X$ is isotropic and $X=G/P$ with $P=P(\mu)$ as above.
It follows that the maximal torus $T$ acts on $G/P$ with fixed points
corresponding to cosets in
\begin{equation}
  \label{eq:not_pi_not}
  \mathcal{Y}=W/W_{\mu}.
\end{equation}

To see this, suppose that $y\in G(k)$ is such that the coset $yP$ is fixed by $T$.   Then
$y^{-1}Ty$ is a maximal torus contained in $P$.   On the other hand, since $P=P(\mu)$ with $\mu:\GG_m\to S\subset T$,
$T$ is also contained in $P$.   Since all maximal tori are conjugate within 
$P$~\cite[Corollary 5.7, p. 496]{SGA3.3}, 
there is a $p\in P(k)$ such that 
\begin{equation}
\label{eq:p_congugator}
p^{-1}Tp=y^{-1}Ty.  
\end{equation} 
Thus $yp^{-1}$ is in the normalizer $N_G(T)$ 
of $T$ and, thus, represents an element $w\in W=W(G,T)$.   If $p'$ is another element of $P$ 
satisfying (\ref{eq:p_congugator}), then $p'p^{-1}$ normalizes $T$.   This implies that 
$p'p^{-1}\in Z_G(\mu)$ by the following.

\begin{lemma}
\label{lemma:parabolic-normalizer}
  For $G$ reductive with maximal torus $T$ and $\mu\in X_*(T)$,
  $P(\mu)\cap N_G(T)\subset Z_G(\mu)$.
\end{lemma}
\begin{proof}
Take $p\in P(\mu)\cap N_G(T)$ and set $\beta(t)=p\mu(t)p^{-1}$.  Since $p$ normalizes $T$,
  $\beta\in X_*(T)$.  Set $a(t)=\mu(t)p\mu (t)^{-1}$.  Since
  $p\in P(\mu )$, $a:=\lim_{t\to 0} a(t)$
  exists.  Therefore
  \begin{equation}
    \label{eq:ExistingLimit}
    \lim_{t\to 0} \beta(t)\mu (t)^{-1} = pa.
  \end{equation}
  In particular, the limit in (\ref{eq:ExistingLimit}) exists.  But,
  since $\beta(t)\mu (t)^{-1}$ is a cocharacter, this is only
  possible if $\beta(t)=\lambda(t)$.
\end{proof}

It follows, therefore, that $y$ and equation (\ref{eq:p_congugator})
determine the class $\pi (y)$ of $yp^{-1}$ in $W/W_{\mu}$.  Moreover, if $y'=yp'$ is another
element of $G$ representing the coset $yP$, then it is easy to check that $\pi (y)=\pi (y')$.

Thus, there is a map $\pi:(G/P)^{T}\to\mathcal{Y}$.
It is not hard to see that the map $\tilde{s}:W\to (G/P)^T$ given by $w\mapsto wP$ induces a map $s:\mathcal{Y}
\to (G/P)^T$ inverse to $\pi$.   Thus $(G/P)^T\cong \mathcal{Y}$.

Now let $yP$ and $y'P$ be two points in the same component $Z$ of $(G/P)^{\lambda}$ which are both
fixed by the $T$ action.  Without loss of generality, we can then assume that $y$ and $y'$ normalize 
$T$.   
By Proposition~\ref{thm:description}, $yP=hy'P$ for some $h\in Z_G(\lambda )$.   
It follows then that $h$ also normalizes $T$ and is, thus, in $W_{\lambda}$.
Thus $\pi$ induces a map  
\begin{equation}
\label{PiNoughtDone}
  q: \pi_0((G/P)^{\lambda})\to \mathcal{X}.
\end{equation}
To see that $q$ is an isomorphism, it suffices to check that $s:\mathcal{Y}\to (G/P)^T$ induces an inverse
map $r:\mathcal{X}\to \pi_0((G/P)^{\lambda}$.  
I leave this verification, which completes the proof of Theorem~\ref{thm:DoubleCosets} to the reader.

\subsection{The Galois action}  
\label{subsec-StarAction}
Let $\Gamma=\Gal(k_{\sep}/k)$ denote the absolute Galois group.  If
$X=G/P(\mu)$ is an isotropic projective homogeneous variety, then $\Gamma$
acts on $G(k_{\sep})$, $T(k_{\sep})$ and $W(k_{\sep})$ stabilizing
$P=P(\mu)$ and, thus, $W_{\mu}$.  It follows that $\Gamma$ acts on the double
coset space $\mathcal{X}$.  Clearly $\Gamma$ also acts on
$(G/P)^{\lambda}(k_{\sep})$, and it is easy to see that the map
$r:\mathcal{X}\to \pi_0(X^{\lambda})$ is an isomorphism of $\Gamma$-sets.

Computing the Galois action on $\pi_0(X^{\lambda})$ can be reduced to
the case of isotropic $X$ using the fact that every reductive
$k$-group has an quasi-split inner form $\Ginn$~\cite[Proposition 16.4.9]{Springer} 
given by a class $\sigma\in H^1(\Gamma, \Gad)$
where $\Gad$ denotes the adjoint group of $G$.  Let $p:G\to\Gad$ be
the canonical quotient map.  It is easy to see that we can arrange
that $T$ is stabilized and $\lambda$ is fixed by $\sigma$.  Then, in
fact, $\sigma$ is in the image of the map 
$$H^1(\Gamma, Z_{\Gad}(p\circ\lambda)\cap N_{\Gad}(p(T)))\to H^1(\Gamma, \Gad).$$
We
have an action of $\Gad$ on $X$ and, thus, a twist $X_{\sigma}$ of $X$
with an action of $\lambda$.  Under the twist, the action of the
Galois group on $W_{\inn}=W(\Ginn, T_{\sigma})$ is given by the
$*$-action~\cite{TitsBoulder}.  Since $G_{\inn}$ is quasi-split,
$X_{\sigma}$ is isotropic and thus corresponds to a parabolic $P(\mu)$
in $\Ginn$.  From the previous section, we then have
$$
\pi_0(X_{\sigma}^{\lambda} ) = W_{\lambda}\backslash W_{\inn} / W_\mu .
$$

\begin{proposition} There is an isomorphism  
$$\pi_0(X^{\lambda})\cong \pi_0(X_{\sigma}^{\lambda})=W_{\lambda}\backslash W_{\inn} / W_\mu $$
of \'{e}tale schemes over $\Spec k$.   (In other words, the above sets are isomorphic as $\Gamma$-sets.)
\end{proposition}  

\begin{proof} It is easy to see that  $\pi_0(X_{\sigma}^{\lambda})$ viewed as an \'{e}tale 
  scheme over $\Spec k$ is the twist of $\pi_0(X^{\lambda})$ by
  $\sigma$.  Since $Z_{\Gad}(p\circ\lambda)$ is geometrically
  connected, it acts trivially on the geometric points of
  $\pi_0(X^{\lambda})$.  Thus the two schemes are isomorphic.
\end{proof}

\section{Explicit Description and Examples}
\label{sec-explicit}

With a little extra work, we can give an explicit description of the
twists and the spaces $Z_i$ appearing in Corollary~\ref{cor:description} 
and Theorem~\ref{thm:RoughAnswer} in terms of the
relevant reflection groups, Dynkin diagrams and root systems.  
This is a generalization of the description appearing in~\cite{CGM2003}.

From the previous section, we know that the quasi-homogeneous schemes
$Z_i$ are in correspondence with the orbits of the $*$-action on the
double cosets in $W_{\lambda}\backslash W_{\inn } / W_{\mu}$.  Over $\ksep$, each such
$Z_i$ decomposes as a disjoint union $Z_w$ over the elements of the
$*$-orbit.  Our goal is then to describe the Tate twist associated to
$Z_w$ and also the projective homogeneous space $Z_w$ in terms of the
Dynkin diagram of $Z(\lambda)$.

In obtaining our description, it will be convenient to consider the
case where $G$ is split first.  Therefore, assume $G$ is split with
maximal torus $T$.  Let $R$ denote the set of roots of $G$.  Choose a
Borel subgroup $B$ or, equivalently, a set $R_+$ of positive roots.
Then $R=R_+\cup R_-$ with $R_-$ the negative roots.  Let $\Sigma$ denote the
corresponding set of simple positive roots.  The Weyl group $W=W(G,T)$
is then generated by the reflections $s_{\alpha}$ in the hyperplanes
defined by the $\alpha\in \Sigma$.  We let $\ell(w)$ denote the corresponding length
function on $W$: $\ell(w)$ is the length $l$ of a minimal expression
$w=s_1s_2\cdots s_l$ of $w$ in terms of the simple roots.

Now let $X$ be a projective homogeneous variety.  We have $X=G/P(\mu)$
for some cocharacter $\mu:\GG_m\to T$ which is non-negative on $R_+$.
(Any cocharacter can be conjugated to a non-negative one.)  Let
$J=\{ \alpha\in \Sigma\, |\, \langle \alpha, \mu\rangle =0\} $.  Then $W_{\mu}$ is
the subgroup of $W$ generated by the $s_{\alpha}$ with $\alpha\in J$.
Accordingly, we will also write $W_J$ for this subgroup.  Now, if $G$
has a non-central cocharacter $\lambda:\GG_m\to T$, there is one which
is non-negative on $\Sigma$.  Thus, setting $I=\{ \alpha\in \Sigma\,
|\, \langle \alpha, \lambda\rangle  =0 \}$, we have $W_{\lambda}\backslash W/
W_{\mu} = W_I\backslash W/ W_J$.  Note that the correspondence
$X\leadsto J$ between isomorphism classes of projective homogeneous
varieties for $G$ and subsets $J\subset \Sigma$ is one-to-one and
onto.  We will call $X$ the homogeneous variety associated to $J$,  and
we will call $J$ the set of roots of $X$.

We now use the result of an exercises  in Humphrey's book on reflection 
groups~\cite[Ex. 1 on p. 20]{HumphreysReflection}.

\begin{exercise}
Any double coset in $W_I\backslash W/ W_J$ has a unique element $b$ of minimal length.  
The element $b$ satisfies the following equivalent properties:
\begin{enumerate}
\item $\ell(bs_{\alpha})=\ell(b)+1$ for $\alpha\in J$ and $\ell(s_{\alpha}b)=\ell(b)+1$ for $\alpha\in I$.
\item $b\alpha>0$ for $\alpha\in J$ and $b^{-1}\alpha >0$ for $\alpha\in I$. 
\end{enumerate}
Moreover, any element $w\in W$ may be written as 
\begin{equation}
  \label{eq:WordForw}
  w=abc
\end{equation}
with $a\in W_I$, $c\in W_J$ and $\ell(w)=\ell(a)+\ell(b)+\ell(c)$.
\end{exercise}

\begin{proof}[Solution]
  The equivalence of the two properties is Lemma~1.6 on p.~12
  of~\cite{HumphreysReflection}.  Let $b$ be an element of minimal
  length in the double coset.  Then $\ell(bs_{\alpha})\geq \ell(b)$ for all $\alpha\in J$,
  and, since lengths either go up by one or down by one upon
  multiplying by a reflection, this implies that $\ell(bs_{\alpha})=\ell(b)+1$
  for all $\alpha\in J$.  Similarly, $\ell(s_{\alpha}b)=\ell(b)+1$ for $\alpha\in I$.  Thus $b$ satisfies both properties (1) and (2). 

Now suppose $w\in W$.  Write $w=abc$ with $\ell(a)$ and $\ell(c)$ minimal.
We can write $a, b$ and $c$ out as reduced words in the simple reflections as follows:
\begin{align*}
   a&=r_1r_2\cdots r_l ,\\
   b&=s_1s_2\cdots s_m ,\\
   c&=t_1t_2\cdots t_n .\\
\end{align*}

Since $a\in W_I$ (resp. $c\in W_J$) and the length function on $W$ restricts to that on $W_I$ (resp. $W_J$), 
we can assume that the $r_i's$ are in $I$ and the $t_i's$ are in $J$. 
So $\ell(w)\leq l+m+n$ and, if $\ell(w)<l+m+n$, there must be a pair of reflections in the word for $w$ that 
can be deleted~\cite[Theorem 1.7]{HumphreysReflection}.   This pair cannot involve one of the $s_i's$ because
otherwise $b$ would not be minimal in the double coset.   On the other hand, it cannot simply involve two of the $r_i's$ because
then our word for $a$ would not be reduced.    Likewise it cannot simply involve two of the $t_i's$.   Finally, the word cannot involve
an $r_i$ and a $t_j$ because we assumed that $\ell(a)$ and $\ell(c)$ were minimal.  

Now suppose $b'$ were another coset representative for $W_IbW_J$ of minimal length.   Then $b'=abc$ for $a$ and $c$ as above.
But, since $\ell(b')=\ell(a)+\ell(b)+\ell(c)$, we must have $a=c=1$.
\end{proof}

\begin{remarks}
  (1) In contrast with the single coset case~\cite[Proposition
  10.10]{HumphreysReflection}, the elements $a$ and $c$ are not
  unique.  This is clear, for example, if $W=\ZZ/2\times \ZZ/2$ and
  $W_I=W_J$ is the first factor of $\ZZ/2$.  
  (2) In the special case
  where $I=J$, this exercise is used in \cite{CGM2003} (see
  Proposition 3.4).  
 (3) What I have stated as the ``exercise'' is
  actually the solution to Humphrey's question which is whether or not
  a minimal element exists.    A more complete version of the exercise is in Bourbaki~\cite[Exercise 1.3, p. 37]{BourbakiCoxeter}. 
\end{remarks}

For a subset $K\subset \Sigma$, let $R_K$ denote the set of roots generated by $K$.    Let $R_K^+=R_K\cap R^+$ (resp.  $R_K^-=R_K\cap R^-$). 
We can now give  our explicit decomposition in the split case.   

\begin{proposition}
  Let $G$ be a split reductive $k$-group with maximal torus $T$ and
  simple roots $\Sigma$.  Let $X=G/P$ be a projective homogeneous variety
  for $G$ with $J$ the corresponding set of simple roots,  let $\lambda$ be
  a non-central cocharacter of $G$ which is non-negative on $\Sigma$ and vanishing precisely on $I\subset \Sigma$ and let
  $E$ be the set of minimal length coset representatives of $W_I\backslash W/
  W_J$ with $W_I$ and $W_J$ as above.  Then, under the
  Bia{\l}ynicki-Birula decomposition for $\lambda$,
  $$
  M(X)= \coprod_{w\in E} M(Z_w)(\ell(w))
  $$
  with $Z_w$ the orbit of $wP$ under $Z(\lambda)$.    The set $I\subset \Sigma$ is the set of simple roots of $Z(\lambda)$.  Moreover, the
roots of $Z_w$ are 
$$
J_w=\{ \alpha\in I\, |\, w^{-1}\alpha\in R_J\}.
$$
\end{proposition}
\begin{proof}
 First note that, if $w\in E$ and $\alpha\in I$, $w^{-1}\alpha\in R^+$.  This is condition (2) of the exercise.  Now the twist $a_w$ associated
to the motive $M(Z(\lambda)wP)$ in the motivic decomposition of Corollary~\ref{cor:description} is the rank of the positive weight space of 
$\lambda$ on $T(G/P)_{wP}$.  This is the same as the rank $r$ of the positive weight space of $\Ad (w^{-1})\lambda$ on $T(G/P)_P$.   
Now $T(G/P)_P$ is naturally identified with $L(G)/ (kR_J^-\oplus L(T)+kR^+)$.   (Here I write $kR$ for the free vector space on the 
set $R$ and view it as a subspace of $L(G)$
in the natural way.)    It follows that $r$ is the number of  negative roots $\alpha$ not in 
$R_J^-$ such that $\langle w^{-1}\lambda, \alpha\rangle  = \langle \lambda, w\alpha\rangle  >0$.    Thus 
\begin{align*}
r&=\#\{ \alpha\in R^--R_J\, |\, w\alpha\in R^+-R_I\} \\
   &=\# (R^+-R_I^+)\cap w(R^--R_J^-)\\
   &=\# (R^+\cap wR^- - R^+\cap wR_J^- - R_I^+\cap wR^-).
\end{align*}
Now $R^+\cap wR_J^-$ and $R_I^+\cap wR^-$ are both empty by part (2) of the exercise.   So
\begin{align*}
r &= \#(R^+\cap wR^-)\\
    &= \ell(w).
\end{align*}

Determining the roots of $Z_w$ is now easy.  We have 
$$Z_w=Z(\lambda)/(Z(\lambda)\cap wPw^{-1}).$$
Set $P_w=Z(\lambda)\cap wPw^{-1}$.  By definition, $P_w$ is a
parabolic subgroup of $Z(\lambda)$.  Note that $P_w$ is actually a {\em
  standard parabolic subgroup} in that it contains the Borel subgroup
$B_{\lambda}=B\cap Z(\lambda)$ where $B$ is the standard Borel subgroup such that
$L(B)=kR^+$.  This follows from the fact that $w^{-1}R_I^+\subset R^+$.  
Now, $L(P)=kR_+\oplus L(T)\oplus kR_J^-$ and $L(Z(\lambda))=kR_I^+\oplus L(T)\oplus kR_I^-$.  
Thus, using the exercise, we have
\begin{align*}
L(P_w)&=(kR_I^+\oplus L(T)\oplus kR_I^-)\cap w(kR_+\oplus L(T)\oplus kR_J^-)\\
              &=kR_I^+ \oplus L(T)\oplus (kR_I^-\cap wkR_J^-)\\
               &=kR_I^+ \oplus L(T)\oplus (kR_I^-\cap wkR_J).
\end{align*}
Now it follows that 
$$J_w=\{ \alpha\in I\, |\, w^{-1}\alpha\in R_J\}.$$
\end{proof}

\subsection{The Tits Index}

To every reductive $k$-group $G$, J.~Tits has associated a
part-algebraic, part-combinatorial object known as the {\em Tits
  index} of $G$.  It consists of the Dynkin diagram for $G$, a graph
with vertices $\Sigma$ corresponding to the simple roots of a given maximal
$k$-torus, together with a Galois action on the vertices preserving a
set $\Sigma_0$ of distinguished vertices.  The subset $\Sigma_0$ consists of the
simple roots orthogonal to a maximal $k$-split torus.  Thus $\Sigma_0$ is the set
of roots of the semi-simple anisotropic kernel.  The group $\Gal(\ksep
/k)$ acts on $\Sigma$ via the $*$-action stabilizing $\Sigma_0$.  When drawing the diagram, the
Galois orbits of the vertices not in $\Sigma_0$ are circled and vertices in
the same $*$-orbit are supposed to be drawn ``close together.''

Now, in this picture, 
projective homogeneous varieties (i.e.,
$k$-defined conjugacy classes of parabolics) are in one-to-one correspondence with
$*$-invariant subset $J$ of $\Sigma$~\cite[(2.5.4)]{TitsBoulder}.
Isotropic projective homogeneous varieties  (i.e., conjugacy classes 
containing a $k$-defined parabolic) 
are in one-to-one correspondence with $*$-invariant subsets
$J$ of $\Sigma$ containing $\Sigma_0$.   
Moreover, it is easy to see that for every $*$-invariant subset $I$ of
$\Sigma$ containing $\Sigma_0$, there is a $k$-defined cocharacter $\lambda$ of $T$
which vanishes on $I$ but is positive on $\Sigma-I$.  In this case, $Z(\lambda)$ is
the Levi component $L_I$ of the parabolic subgroup $P_I$ associated to $I$.
(The necessary argument is given in the proof of \cite[Lemma
15.1.2]{Springer}).

 We now have the following result which follows from the previous proposition by standard methods
of descent.

\begin{theorem}
\label{thm-Final}
Suppose $I$ and $J$ are $*$-invariant subsets of $\Sigma$ with $I\supset \Sigma_0$.    Let $E$ be the set of minimal length coset representatives 
for $W_I\backslash  W_{\inn} / W_J$, and let $\overline{E}$ be the 
set of orbits of $E$ under the $*$-action.   Let $X$ be the projective homogeneous variety associated to $J$.    We have
$$
M(X)=\coprod_{\overline{w}\in \overline{E}}  M(Z_{\overline{w}})(\ell(w))
$$
where $Z_{\overline{w}}$ is a projective quasi-homogeneous scheme for the reductive group $L_I $ $(=Z(\lambda))$.  
Moreover, the base change of $Z_{\overline{w}}$ to $\ksep$ is a disjoint union
$$
Z_{\overline{w}} \otimes \ksep = \coprod_{w\in \overline{w}} Z_w 
$$
where $Z_w$ is the projective homogeneous variety for $L_I\otimes \ksep$ corresponding to the subset 
$$
J_w=\{ \alpha\in I\, |\, w^{-1}\alpha\in R_J\}.
$$
\end{theorem}

\begin{remark}
If the $*$-orbit $\overline{w}$ consists of one element $w$, then $Z_{\overline{w}}$
is the projective homogeneous variety corresponding to $J_w$.  In particular, its structure as a 
$k$-variety is determined by the combinatorics.   At any rate, since $Z_{\overline{w}}$ is a projective quasi-homogeneous
variety, it is a subvariety of the variety of parabolics of the reductive group $L_I$.
\end{remark}

\begin{example}
\label{example-TwoDIsotropic}
  We now return to the example of (\ref{rem-InterestingExample}).
  Here we have $G=PSO(q)$ with $q=x_1^2+\cdots +x_{2n-2}^2 + yz$  (with $k=\RR$) and $X$ the
  projective quasi-homogeneous scheme of two-dimensional subspaces.
  As long as $n\geq 3$, $X$ is a
  projective homogeneous variety for $G$.  (When $n=2$, $X$ has two geometric
  components.)

 Now the Dynkin diagram of $G$ (decorated as in the Tits index) is the following picture.

 $$
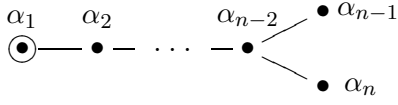
\xy <1cm,0cm>:
(0,0)*+\cir<5pt>{};
(0,.4)*+{\alpha_1}; (1, .4)*+{\alpha_2}; 
(0,0)*+{\bullet };  (1,0)*+{\bullet} **@{-};
(1.6,0)*+{} **@{-};
(1.8,0)*+{.}; (2,0)*+{.}; (2.2,0)*+{.};
(2.4,0)*+{}; (3,0)*+{\bullet} **@{-};
(3,.4)*+{\alpha_{n-2}};
(3.2, .1);   (4,.5)*+{\bullet} **@{-};
(3.2, -0.1); (4,-.5)*+{\bullet} **@{-};
(4.6, .5)*+{\alpha_{n-1}};
(4.5, -.5)*+{\alpha_n};
\endxy
$$

The $*$-action exchanges $\alpha_{n-1}$ with $\alpha_n$ and leaves all other
roots fixed.  When $n\geq 4$, the set $J\subset\Sigma$ corresponding to $X$ is
$\Sigma-\{\alpha_2\}$.  For $n=3$, $J=\Sigma-\{\alpha_2,\alpha_3\}$.  We assume $n\geq 4$ at first and
sketch the case where $n=3$ (where the $*$-action plays a significant
role) at the end of this example.

If we set $I=\Sigma_0=\Sigma-\{\alpha_1\}$, then we are in the setting of
Theorem~\ref{thm-Final}.  Write $s_i=s_{\alpha_i}$ for the generators of
the Weyl group.  From the theory of Coxeter complexes (\cite[\S{}
1.15]{HumphreysReflection}), we see that $W_I\backslash W$ is identified with
the set of vectors of the form $\pm e_i$ in the real vector space
$\RR^n=\RR e_1+\cdots +\RR e_n$.  The action of $W$ on the right is given
by $e_is_j=e_{(j,j+1)i}$ for $1\leq j<n$ where $(j,j+1)$ denotes the
transposition in the symmetric group exchanging $j$ and $j+1$.  For
$j=n$, we have $e_is_n=-e_{(n-1,n)i}$.  Now it is fairly easy to see
that the cosets in $W_I\backslash W$ containing elements of $E$ are $e_1=W_I1$,
$e_3=W_Is_1s_2$, and $-e_2=W_Is_1s_2\cdots s_{n-2} s_{n-1}s_ns_{n-2}\cdots
s_3s_2$.  It is also easy to see that the representatives listed are
in fact the elements of $E$.  That is,
$$
E=\{1, s_1s_2, s_1s_2\cdots
s_{n-2} s_{n-1}s_ns_{n-2}\cdots s_2 
\}.
$$ 
Writing $w_1, w_2$ and $w_3$ for the elements listed in order, we have
$$
\ell(w_1)=0, \ell(w_2)=2, \ell(w_3)=2n-3.
$$  

Clearly $J_{w_1}=I\cap
J=\{\alpha_3, \ldots, \alpha_n\}$.  To compute $J_{w_2}$, note that
\begin{align*}
w_2^{-1}\alpha_2 &=s_2s_1\alpha_2\\
                     &=s_2(\alpha_1+\alpha_2)\\
                     &=(\alpha_1+\alpha_2)-\alpha_2\\
                     &= \alpha_1\in J.
\end{align*}
Thus $\alpha_2\in J_{w_2}$.  A similar computation shows that, for $i>2$,
$w_2^{-1}\alpha_i\in J_{w_2}$ if and only if $\alpha_i$ is not connected to $\alpha_2$
in the Dynkin diagram.  Thus, for $n\geq 5$, $J_{w_2}=\{ \alpha_2, \alpha_4,\ldots,
\alpha_n\}$.  However, for $n=4$, $J_{w_2}=\{ \alpha_2\}$.  

Finally, to compute $J_{w_3}$, note that
$w_3^{-1}\alpha_i=\alpha_i$ for $2<i\leq n-2$, $w_3^{-1}\alpha_2=\alpha_1+\alpha_2$,
$w_3^{-1}\alpha_{n-1}=\alpha_n$ and $w_3^{-1}\alpha_n=\alpha_{n-1}$.  
It follows that $J_{w_3}=J_{w_1}=I-\{\alpha_2\}$. 

Putting all of this together, we have the decomposition
\begin{equation}
  \label{eq:FinalExample}
  M(X)=M(Q)\oplus M(Y)(2)\oplus M(Q)(2n-3)
\end{equation}
where $Q$ is the motive of a quadric  of isotropic lines for the quadratic form
$q'=x_1^2 +\cdots +x_{2n}^2$, and $Y$ is isomorphic to the space of
isotropic planes for $q'$.

When $n=3$, $W_I\backslash W$ (viewed in terms of the Coxeter complex) has four $W_J$ orbits containing the 
minimal $W_I$-cosets
\begin{align*}
  e_1&=W_Iw_1,\quad w_1=1;\\
  e_3&= W_Iw_2,\quad w_2=s_1s_2;\\
  -e_3&=W_Iw_3,\quad w_3=s_1s_3;\\
  -e_2&= W_Iw_r,\quad w_4=s_1s_2s_3.
\end{align*}

The elements $w_i$ are written in reduced form so we have $\ell(w_1)=1, \ell(w_2)=\ell(w_3)=2$ and $\ell(w_4)=3$.
Note that $w_2$ and $w_3$ are conjugate under the $*$-action.
Clearly $J_{w_1}=\emptyset$ and some computation shows that $J_{w_4}=\emptyset$ as well.    
On the other hand, $J_{w_2}= \alpha_2$ while $J_{w_3}=\alpha_3$.

It turns out then that $Z_{\overline{w_1}}=Z_{\overline{w_4}}=Q$ and $Z_{\overline{w_3}}=Y$ 
with $Q$ and $Y$ as in~\ref{eq:FinalExample}.   In other words, we obtain the same decomposition 
as in the case $n=4$.   However,  we learn that $Y\otimes \CC=Z_{w_2}\coprod Z_{w_3}$ with  $Z_{w_2}=Z_{w_3}=\PP^1$.

It is also possible (and perhaps easier) to work out (\ref{eq:FinalExample}) directly using the geometry
of the $\GG_m$-action on $X$ and the weight decomposition of $TX$ at the various fixed loci as suggested
in Exercise~\ref{rem-InterestingExample}.
\end{example}

\section{Acknowledgments}   
It is a pleasure to thank Chernousov, Gille, and Merkurjev as well as
P.~Belkale, N.~Fakhruddin and O.~Gabber for extremely helpful
discussions and H.~Furdyna and B.~Brosnan for pointing out several
typos.  I would also like to thank B.~Shipman for sending me a reprint of 
her paper~\cite{ShipmanFixedPoints} describing the connected
components of $(G/B)^T$ which helped me greatly to guess
Theorem~\ref{thm:DoubleCosets}, J.~Carrell for sending me useful
references to his paper ``Torus actions and cohomology'' in
\cite{BbCarrell} and B.~K\"{o}ck for telling me about del Ba\~{n}o's
work.



\begin{thebibliography}{10}

\bibitem{SGA3.3}
{\em Sch\'emas en groupes. {III}: {S}tructure des sch\'emas en groupes
  r\'eductifs}.
\newblock S\'eminaire de G\'eom\'etrie Alg\'ebrique du Bois Marie 1962/64 (SGA
  3). Dirig\'e par M. Demazure et A. Grothendieck. Lecture Notes in
  Mathematics, Vol. 153. Springer-Verlag, Berlin, 1962/1964.

\bibitem{BialynickiBirulaDecomposition}
A.~Bia{\l}ynicki-Birula.
\newblock Some theorems on actions of algebraic groups.
\newblock {\em Ann. of Math. (2)}, 98:480--497, 1973.

\bibitem{BialynickiBirulaFiltration}
A.~Bia{\l}ynicki-Birula.
\newblock Some properties of the decompositions of algebraic varieties
  determined by actions of a torus.
\newblock {\em Bull. Acad. Polon. Sci. S\'er. Sci. Math. Astronom. Phys.},
  24(9):667--674, 1976.

\bibitem{BbCarrell}
A.~Bia{\l}ynicki-Birula, J.~B. Carrell, and W.~M. McGovern.
\newblock {\em Algebraic quotients. {T}orus actions and cohomology. {T}he
  adjoint representation and the adjoint action}, volume 131 of {\em
  Encyclopaedia of Mathematical Sciences}.
\newblock Springer-Verlag, Berlin, 2002.
\newblock Invariant Theory and Algebraic Transformation Groups, II.

\bibitem{borel-tits}
Armand Borel and Jacques Tits.
\newblock Groupes r\'eductifs.
\newblock {\em Inst. Hautes \'Etudes Sci. Publ. Math.}, (27):55--150, 1965.

\bibitem{BourbakiCoxeter}
N.~Bourbaki.
\newblock {\em \'{E}l\'ements de math\'ematique. {F}asc. {XXXIV}. {G}roupes et
  alg\`ebres de {L}ie. {C}hapitre {IV}: {G}roupes de {C}oxeter et syst\`emes de
  {T}its. {C}hapitre {V}: {G}roupes engendr\'es par des r\'eflexions.
  {C}hapitre {VI}: syst\`emes de racines}.
\newblock Actualit\'es Scientifiques et Industrielles, No. 1337. Hermann,
  Paris, 1968.

\bibitem{brosnan-documenta}
Patrick Brosnan.
\newblock A short proof of rost nilpotence via refined correspondences.
\newblock {\em Doc. Math.}, 8:79--96 (electronic), 2003.

\bibitem{CarrellSommese}
James~B. Carrell and Andrew~John Sommese.
\newblock Filtrations of meromorphic {${\bf C}\sp{\ast} $} actions on complex
  manifolds.
\newblock {\em Math. Scand.}, 53(1):25--31, 1983.

\bibitem{CGM2003}
Vladimir Chernousov, Stefan Gille, and Alexander Merkurjev.
\newblock
  \href{http://www.math.ucla.edu/~merkurev/papers/nilpotence9.dvi}{Motivic
  decomposition of isotropic projective homogeneous varieties}.
\newblock Currently available at
  \url{http://www.math.ucla.edu/~merkurev/papers/nilpotence9.dvi}.
\newblock To appear in Duke Math. Journal.

\bibitem{delBano}
Sebastian del Ba{\~n}o.
\newblock On the {C}how motive of some moduli spaces.
\newblock {\em J. Reine Angew. Math.}, 532:105--132, 2001.

\bibitem{fulton}
William Fulton.
\newblock {\em Intersection theory}.
\newblock Springer-Verlag, Berlin, second edition, 1998.

\bibitem{HartshorneAG}
Robin Hartshorne.
\newblock {\em Algebraic geometry}.
\newblock Springer-Verlag, New York, 1977.
\newblock Graduate Texts in Mathematics, No. 52.

\bibitem{HesselinkConcentration}
Wim~H. Hesselink.
\newblock Concentration under actions of algebraic groups.
\newblock In {\em Paul Dubreil and Marie-Paule Malliavin Algebra Seminar, 33rd
  Year (Paris, 1980)}, volume 867 of {\em Lecture Notes in Math.}, pages
  55--89. Springer, Berlin, 1981.

\bibitem{HumphreysReflection}
James~E. Humphreys.
\newblock {\em Reflection groups and {C}oxeter groups}, volume~29 of {\em
  Cambridge Studies in Advanced Mathematics}.
\newblock Cambridge University Press, Cambridge, 1990.

\bibitem{IversenFixed}
Birger Iversen.
\newblock A fixed point formula for action of tori on algebraic varieties.
\newblock {\em Invent. Math.}, 16:229--236, 1972.

\bibitem{Jannsen}
Uwe Jannsen.
\newblock Motives, numerical equivalence, and semi-simplicity.
\newblock {\em Invent. Math.}, 107(3):447--452, 1992.

\bibitem{karpenko-flag}
N.~A. Karpenko.
\newblock Cohomology of relative cellular spaces and of isotropic flag
  varieties.
\newblock {\em Algebra i Analiz}, 12(1):3--69, 2000.

\bibitem{Knus}
Max-Albert Knus.
\newblock {\em Quadratic and {H}ermitian forms over rings}, volume 294 of {\em
  Grundlehren der Mathematischen Wissenschaften [Fundamental Principles of
  Mathematical Sciences]}.
\newblock Springer-Verlag, Berlin, 1991.
\newblock With a foreword by I. Bertuccioni.

\bibitem{Kock}
Bernhard K{\"o}ck.
\newblock Chow motif and higher {C}how theory of {$G/P$}.
\newblock {\em Manuscripta Math.}, 70(4):363--372, 1991.

\bibitem{Lam}
T.~Y. Lam.
\newblock {\em The algebraic theory of quadratic forms}.
\newblock W. A. Benjamin, Inc., Reading, Mass., 1973.
\newblock Mathematics Lecture Note Series.

\bibitem{Manin}
Ju.~I. Manin.
\newblock Correspondences, motifs and monoidal transformations.
\newblock {\em Mat. Sb. (N.S.)}, 77 (119):475--507, 1968.

\bibitem{GIT}
D.~Mumford, J.~Fogarty, and F.~Kirwan.
\newblock {\em Geometric invariant theory}, volume~34 of {\em Ergebnisse der
  Mathematik und ihrer Grenzgebiete (2) [Results in Mathematics and Related
  Areas (2)]}.
\newblock Springer-Verlag, Berlin, third edition, 1994.

\bibitem{RostPfister}
Markus Rost.
\newblock
  \href{http://www.mathematik.uni-bielefeld.de/~rost/data/motive.pdf}{The
  motive of a Pfister form}.
\newblock Currently available at
  \url{http://www.mathematik.uni-bielefeld.de/~rost/data/motive.pdf}.

\bibitem{ShipmanFixedPoints}
Barbara~A. Shipman.
\newblock On the fixed-point sets of torus actions on flag manifolds.
\newblock {\em J. Algebra Appl.}, 1(3):255--265, 2002.

\bibitem{Springer}
T.~A. Springer.
\newblock {\em Linear algebraic groups}, volume~9 of {\em Progress in
  Mathematics}.
\newblock Birkh\"auser Boston Inc., Boston, MA, second edition, 1998.

\bibitem{TitsBoulder}
J.~Tits.
\newblock Classification of algebraic semisimple groups.
\newblock In {\em Algebraic Groups and Discontinuous Subgroups (Proc. Sympos.
  Pure Math., Boulder, Colo., 1965)}, pages 33--62. Amer. Math. Soc.,
  Providence, R.I., 1966, 1966.

\end{thebibliography}
\def\noopsort#1{} \def\cprime{$'$} \def\noopsort#1{} \def\cprime{$'$}

\end{document}